\documentstyle{article}
\newtheorem{thm}{Theorem}
\newtheorem{prop}[thm]{Proposition}
\newtheorem{lem}[thm]{Lemma}
\newtheorem{cor}[thm]{Corollary}
\newtheorem{main}{Theorem}

\newtheorem{defn}{Definition}

\newenvironment{xpl}{\medskip 
\noindent {\bf  Example.}}{\mbox{}$\quad$\bigskip}
\newenvironment{proof}{\medskip 
\noindent {\bf Proof.}}{\hfill \rule{.5em}{1em}\mbox{}\bigskip}
\def\bea{\begin{eqnarray*}}
\def\eea{\end{eqnarray*}}
\def\bel{\begin{eqnarray}}
\def\eel{\end{eqnarray}}
\def\be{\begin{equation}}
\def\ee{\end{equation}}

\def\Bbb{\bf }
\title{Yamabe Invariants and Spin$^c$ Structures}
\author{
Matthew J. Gursky\thanks{Supported 
in part by  NSF grant DMS-9623048}\\Indiana University  
\\and\\ 
Claude LeBrun
\thanks{Supported 
in part by  NSF grant DMS-9505744}\\SUNY Stony Brook}

\begin{document}
\maketitle

\begin{abstract} 
The {\em Yamabe invariant} of a smooth compact manifold
is by definition the supremum
 of the  scalar curvatures of 
unit-volume Yamabe metrics on the manifold. 
For  an explicit infinite class of 
$4$-manifolds, we show  that  
this invariant is positive but strictly less than
that of the  $4$-sphere. This is done by
using  spin$^c$ 
Dirac operators to control
the lowest eigenvalue of a
perturbation of the Yamabe Laplacian. 
These results dovetail perfectly with those derived from
 the perturbed 
Seiberg-Witten equations \cite{leb4},
but  the present method is much  more elementary
in spirit. 
 \end{abstract}
 
 \section{Introduction}
There is a natural diffeomorphism invariant
\cite{okob,sch1} 
which arises
from a variational problem for 
  the  total scalar curvature of Riemannian metrics on 
any given  compact smooth $n$-manifold $M$.
Observe that the group of smooth positive 
functions  $u: M\to {\Bbb R}^+$
acts on the space of smooth Riemannian metrics
$g$ by conformal rescaling $g\mapsto u^2g$.
The  {\em conformal class} of a 
Riemannian metric $g$ is by definition 
 the orbit
 $\gamma = [g]$ of this action which 
  contains $g$.   
Let ${\cal C}(M) =\{ \gamma \}$ denote the set of
conformal classes of metrics 
on $M$. We may then define an invariant of the 
smooth compact manifold $M$ by 
setting
$$Y(M) := 
\sup_{\gamma \in {\cal C}(M)}\inf_{g\in \gamma} \frac{\int_M
s_g~d\mu_g}{\left(\int_M 
d\mu_g\right)^{\frac{n-2}{n}}},
$$
where $s_g$ and $d\mu_g$ respectively
 denote the scalar curvature and volume measure of
the Riemannian metric $g$.
We will call this the {\em Yamabe invariant}
of $M$.

To put this definition in context, recall \cite{bes}
that, for $n > 2$, a Riemannian metric is Einstein
iff it is a critical point of the functional 
$$g \mapsto  {\cal S}(g) :=
\frac{\int_M 
s_g~d\mu_g}{\left(\int_M 
d\mu_g\right)^{\frac{n-2}{n}}}.$$ 
The functional $\cal S$, however, is neither bounded above nor below,
so one cannot hope to find a critical point by
either minimizing or maximizing it. However, 
the restriction of $\cal S$ to any conformal
class {\sl is} bounded below, and a 
remarkable theorem
\cite{aubin,lp,rick} 
of Yamabe, Trudinger, Aubin, and Schoen 
asserts that each conformal class $\gamma$
contains   metrics $g$, called {\em Yamabe metrics},
which attain the minimum value
$$Y_{\gamma } = \inf_{g\in \gamma  } {\cal S}(g).$$
This number is called the {\em Yamabe constant} of the
conformal class $\gamma$. 
A simple and beautiful argument of Aubin \cite{aubin}
shows  that  
$Y_\gamma\leq Y(S^n)= n(n-1) V^{2/n}_n$
for any conformal class on any $n$-manifold,
 where $V_n$ is the 
volume of the standard  metric on  $S^n$.  
Thus the
scalar curvatures of unit-volume Yamabe
metrics on $M$ are bounded above, and their
least upper bound is a real number 
  $Y(M)\leq Y(S^n)$. Of course, this by no means 
guarantees\footnote{
A plausible folk
tradition   maintains that 
 it was  Yamabe himself 
who first  considered the question of whether 
$Y(M)$ can be realized as the scalar curvature of 
a unit-volume Yamabe metric. 
 However, $Y(M)$ seems to 
make its first  published appearance in 
an important 
 paper of   O. Kobayashi \cite{okob}, who called it
the mu invariant. Other authors   \cite{and,sch1}
have elsewhere called it the  sigma  constant.}
 that $Y(M)$ is a critical 
value of 
$\cal S$. 
Indeed, there are
 many low-dimensional examples \cite{bes} 
of manifolds which do not admit any  Einstein metrics.

A conformal class $\gamma$ contains a metric
$g\in \gamma$ with positive scalar curvature 
iff $Y_\gamma > 0$, so the Yamabe invariant $Y(M)$ is
positive iff $M$ carries a metric of positive 
scalar curvature. Now there is a substantial 
body of results \cite{lic,hit,GL,SY,LM,stolz} 
concerning manifolds which admit metrics of 
positive scalar curvature, and  these
results may be understood   as simply providing 
one kind of estimate for  
Yamabe invariants. The bulk  of this literature
consists of variations on a theme 
 of Lichnerowicz: 
on a spin manifold of positive scalar curvature,
the Dirac operator must have  index zero.
In the present article, we will show that essentially
the same method, applied to the {\em twisted}  Dirac operators
associated with spin$^c$ structures, can be
used   to calculate the Yamabe invariant
for  many 4-manifolds.

For rather mysterious  reasons, 
the Yamabe invariant seems to be most easily  computed  
in low dimensions. In dimension 2, for example, 
the  Gauss-Bonnet theorem asserts that $Y(M)$ is 
  a multiple of the Euler characteristic.
In dimension 3, Anderson \cite{and} has announced  a
computation of the 
  Yamabe invariants of all
irreducible 3-manifolds with infinite $\pi_1$.
And in  dimension 4, which will be our field of concern, 
  the advent of 
Seiberg-Witten theory 
  \cite{witten} has made it possible \cite{leb3,leb5} to compute 
the Yamabe invariants of most complex algebraic surfaces.
One  remarkable feature that emerges is that 
  $Y(M)$ often  distinguishes between 
 different smooth structures on the same topological
4-manifold. 

The Seiberg-Witten method, however, is most finely tuned  
to manifolds with $Y(M)\leq 0$, although a perturbed version
 can be used \cite{leb4} to show, for example, 
 that $Y({\Bbb CP}_2)=12\sqrt{2}\pi
< 8\sqrt{6}\pi =Y(S^4)$. In this article,  
 the last result will be reproved by a much simpler method,  
while at the same time proving the following 
substantial generalization:

\begin{main}
The Yamabe invariant of 
 ${\Bbb CP}_2$ is  unaltered by $0$-surgeries:  
$$Y({\Bbb CP}_2\# m (S^1\times S^3))=Y( {\Bbb CP}_2) = 12\sqrt{2}\pi $$
for all $m \geq 0$. In particular, these 
projective planes with handles  all 
have Yamabe invariant strictly less than 
$Y(S^4)=8\sqrt{6}\pi$.
\end{main} 

The same reasoning also proves the following:

\begin{main} Let $k\in \{1,2,3\}$, and let 
$m$ be any natural number. Then 
$$12\sqrt{2}\pi \leq 
Y(k{\Bbb CP}_2\#  m (S^1\times S^3)) \leq 4\pi \sqrt{2k+16}.$$
In particular, these connected sums of ${\Bbb CP}_2$'s and
$S^1\times S^3$'s  all 
have Yamabe invariant strictly less than 
$Y(S^4)$.
\end{main} 

For K\"ahler-type complex surfaces of 
Kodaira dimension $\geq 0$, 
  Seiberg-Witten theory  allows one
to show  \cite{leb5} that   $Y(M)$ is 
a bimeromorphic invariant --- i.e. it is unchanged by 
blowing up and down. One might  therefore blithely 
suppose that the same is true of  {\em all}  
complex surfaces. However, the present methods show that this
 supposition simply does not hold water:

\begin{main} The Hopf surface   and its
one-point blow-up 
have different Yamabe invariants. Thus the Yamabe
invariant is {\em not} a bimeromorphic invariant for 
complex surfaces of    class $VII$. 
\end{main}

The key  inequality  for the Yamabe
invariant developed here is sensitive only
to  homotopy type 
rather than to  diffeomorphism type.
The present methods   are thus 
oblivious  to the 
deeper aspects of 4-manifold topology
detected by  Seiberg-Witten invariants. 
Nonetheless,   several peculiarities of dimension $4$ ---
e.g.  the conformal invariance of
harmonic 2-forms --- will play a   a crucial r\^ole.
It   thus remains to be seen whether 
the ideas developed in this article have 
any ramifications 
in higher dimensions.

\section{Perturbed Yamabe Laplacians}

Let $(M^4,g)$ be a smooth 4-dimensional Riemannian manifold.
The   {\em  Yamabe Laplacian}  of
$g$ will mean the elliptic operator
\begin{equation}
\Box_g= 6\Delta_g+ s_g \label{one}
\end{equation}
acting on functions on $M$.
Here $s_g$ denotes the scalar curvature of $g$,
and  $\Delta =d^*d=-\mbox{div } \mbox{grad}$ 
is the    (positive) Laplace-Beltrami operator 
  of  $g$. Considered as a map between spaces of 
functions (or densities) of appropriate conformal 
weight, 
this operator is {\em conformally invariant}; namely,
if $\tilde{g}=u^2g$ for some positive $C^2$ function $u$, then 
\begin{equation}
\Box_{\tilde{g}}\varphi= u^{-3}\Box_g (u\varphi)
\label{two}
\end{equation}
for any function $\varphi$.
The   geometric essence of this
statement is the fact  that 
the  scalar curvature transforms under conformal rescalings
according
to the rule
$$s_{\tilde{g}}= u^{-3}\Box_g u.$$
 
 Let $E\subset \otimes^2T^*M$ be the real line  bundle 
spanned by the metric $g$. Evidently,
this depends only on the conformal class
$\gamma $, and conversely the conformal class is
uniquely determined by $E$. 
A section $f$ of $E$ may simply be thought of 
as a real valued  function on $M$ which transforms
according to the rule
$$f \mapsto \tilde{f}= u^{-2} f$$
when the metric $g$ is replaced by $\tilde{g}=u^2g$,
since this transformation rule ensures that
$\tilde{f}\tilde{g}=fg$.  Sections of $E$ will
therefore be called {\em functions of conformal weight
$-2$.}

\begin{xpl}
Let $\omega$ be a smooth 2-form. The function 
$$f= |\omega |_g$$
then transforms
according to the rule
$$f\mapsto \tilde{f}=u^{-2}f$$
when $g \mapsto \tilde{g}=u^2g$. Thus
$f$ is a function of conformal weight $-2$. Notice
that while $f^2$ is smooth, $f$    will typically 
merely be  Lipschitz if the 
locus where $\omega$ vanishes is non-empty. 
\end{xpl}

\begin{lem}
Let $\gamma $ be a smooth conformal class on a 4-manifold $M$, and let
$f$ be a function of conformal weight $-2$ on $M$. Then
the operator 
$\diamondsuit_{g}=\Box_g - f$
transforms according to the rule
$$\diamondsuit_{\tilde{g}}\varphi= u^{-3}\diamondsuit_g(u\varphi)$$
when $g$ is replaced by $\tilde{g}=u^2g.$
\end{lem}
\begin{proof}
We have 
\bea
\diamondsuit_{\tilde{g}}\varphi&=& \Box_{\tilde{g}}\varphi - 
\tilde{f} \varphi
\\&=& u^{-3} \Box_g(u\varphi) -u^{-2} f \varphi\\&=&
 u^{-3} [\Box_g(u\varphi) -u  f \varphi ]\\&=&u^{-3}\diamondsuit_g(u\varphi)
\eea
by the conformal invariance of the 
Yamabe Laplacian.
\end{proof}

\begin{defn}
Let $g$ be a metric on $M$, and let $f$ be a function
of conformal weight $-2$. The 
{\em modified scalar curvature} of the pair $(g,f)$
will mean the function $\sigma = \sigma_{(g,f)}= s - f$,
where $s=s_g$ is the usual scalar curvature of $g$. 
\end{defn}

\begin{lem}
Under conformal changes $g\mapsto \tilde{g}=u^2g$,
the modified scalar curvature transforms
according to the rule $\sigma \mapsto \tilde{\sigma}=
u^{-3}\diamondsuit_gu.$
\end{lem}
\begin{proof}
Indeed, $\sigma_{(g,f)} = \diamondsuit_g(1).$ By the  
previous lemma, we therefore have $\sigma_{(\tilde{g},\tilde{f})}=
\diamondsuit_{\tilde{g}}(1)= u^{-3}\diamondsuit_gu.$
\end{proof}

\begin{prop}
Let $g$ be a smooth Riemannian metric  on a compact
smooth 4-manifold $M$,
and let $f\in C^{0,\alpha}(M,E)$,
$\alpha \in (0,1)$,  be a   H\"older continuous
function of conformal weight $-2$. 
 Then  there is a 
conformally related metric $\tilde{g}=u^2g$ of class  $C^{2,\alpha}$
whose modified   scalar curvature   satisfies
$\tilde{\sigma} > 0$, $\tilde{\sigma} < 0$, or 
$\tilde{\sigma} \equiv 0.$
 Moreover, these three possibilities
are mutually exclusive. \label{try} 
\end{prop}
\begin{proof}
Let $\lambda_g$ be the lowest eigenvalue of $\diamondsuit_g$:
$$\lambda_g=\inf_{\begin{array}{c}u \in L^2_1\\ 
\|u\|_{L^2}=1\end{array}}
\langle \diamondsuit_gu , u \rangle_{L^2(g)}.$$
Let $u$ be a non-zero
 eigenfunction of $\diamondsuit_g$ corresponding to this eigenvalue:
$$\diamondsuit_gu=\lambda_g u.$$
By the interior
Schauder estimates \cite[p.109]{gt}, 
$u$ is of class $C^{2,\alpha}$.
By the minimum principle \cite[p.35]{gt}, $u\neq 0$, 
so 
  $\tilde{g}=u^2g$
is a $C^{2,\alpha}$  metric  
conformal to $g$. Its modified scalar curvature
is 
$$\tilde{\sigma} = u^{-3} \diamondsuit_gu = u^{-2}\lambda_g,$$
and so is strictly positive, strictly negative, or identically zero,
exactly as promised. 

Now notice that  $\diamondsuit_{\tilde{g}}= 6\Delta_{\tilde{g}}+
\tilde{\sigma}$. Thus, for any 
positive $C^2$ function  $v$,   the 
  modified scalar curvature  
of   $v^2\tilde{g}$ is at most 
$v^{-2}\tilde{\sigma}$ at the  minima of $v$, and  
at least $v^{-2}\tilde{\sigma}$
at the    maxima of $v$.  The   three possibilities 
under discussion 
are therefore  mutually exclusive. 
\end{proof}

Notice that the $L^2$ norm 
$$
\| \omega\|_{2} = \left( \int_M |\omega |_g^2d\mu_g
\right)^{1/2}
$$
of a 2-form $\omega$ on any compact 4-manifold $M$ is 
{\em conformally invariant}; that is, it depends only on 
the conformal class $\gamma = [g]$ of the metric. 

\begin{cor}\label{yup}
Let $\gamma $ be a smooth conformal class on  
  a smooth compact   4-manifold $M$, and  
let $\omega$ be a differentiable 2-form on $M$.
Then   one of the following must hold: 
\begin{itemize}
\item  there is a $C^\infty$  metric ${g}\in \gamma $ 
with scalar curvature ${s} > |\omega |_{{g}}$;
or 
\item  
$Y_\gamma <  \| \omega\|_{2}$; or 
\item  
$Y_\gamma = \| \omega\|_{2}$,
and there is a ($C^\infty$) Yamabe   metric $g\in \gamma$
 with $s = |\omega | \equiv \mbox{\rm const}$. 
In particular, this happens only if 
$\omega$ is   nowhere zero or vanishes identically.
\end{itemize}
\end{cor}
\begin{proof} Let $f=|\omega|$. The corresponding 
modified scalar curvature $\sigma = s - |\omega |$
then defines a continuous map from the 
Banach space of $C^2$ metrics in $\gamma$ to the 
Banach space of $C^0$ functions. Thus the set of $C^2$ metrics
in $\gamma$ with $\sigma > 0$ is therefore $C^2$ open. However, the
smooth metrics in $\gamma$  are dense in the $C^2$ metrics. Thus,
if there is no smooth metric in $\gamma$ with 
$s > |\omega |$, there cannot be a $C^{2,\alpha}$ 
metric with $s > |\omega |$ either. But by Proposition \ref{try},
this happens precisely if there is instead  a 
$C^{2,\alpha}$ metric $g\in \gamma$ with $s \leq |\omega |$.

If 
 the latter happens, we then have a metric $g\in \gamma$
for which 
$$\frac{\int s~d\mu}{\sqrt{\int d\mu}}
\leq \frac{\int |\omega | d\mu}{\sqrt{\int d\mu}}
 \leq \sqrt{\int |\omega |^2 d\mu} =\|\omega \|_2, $$
so that the definition of the Yamabe constant yields
$$Y_\gamma \leq  \|\omega\|_2. $$
If equality holds, moreover, the  metric $g$
is a Yamabe metric, and satisfies $s\equiv |\omega|$.
Since $g$ is a Yamabe metric, it is smooth   and has 
constant scalar curvature. In particular,
$|\omega |\equiv s$ must be constant. 
\end{proof}

\section{Polarizations}

Let $M$ be a compact oriented 4-manifold, and let
$\gamma$ be a conformal class on $M$. Then
the orientation and conformal structure  induce 
a Hodge star operator 
$$\star: \Lambda^2 \to \Lambda^2$$
on the bundle of 2-forms.
That is to say, the Hodge star operator 
on middle-dimensional forms determined
by any metric $g\in \gamma$ actually depends only 
on the conformal class $\gamma$.
This linear endomorphism of $\Lambda^2$ 
satisfies $\star^2=1$, so that
we have an eigenspace decomposition
$$\Lambda^2=\Lambda^+\oplus \Lambda^-$$
depending only on $\gamma$ and the orientation.
The factors $\Lambda^\pm$, corresponding to the 
eigenvalues $\pm 1$,  are vector bundles of
rank 3, and reversing the orientation of 
$M$ just interchanges them. Sections of 
$\Lambda^+$ are called self-dual 2-forms,
whereas sections of $\Lambda^-$ are called
anti-self-dual 2-forms. 

Now the Hodge theorem tells us that 
$$H^2(M,{\Bbb R}) =\{ \varphi \in \Gamma (\Lambda^2) ~|~
d\varphi = 0, ~ d\star \varphi =0 \} .$$
Since $\star$ defines an involution of 
the right-hand side, however, 
we therefore have a direct sum decomposition
$$H^2(M, {\Bbb R}) = {\cal H}^+\oplus {\cal H}^-,$$
where
$${\cal H}^\pm= \{ \varphi \in \Gamma (\Lambda^\pm) ~|~
d\varphi = 0\} $$
are the spaces of self-dual and anti-self-dual harmonic forms.
Given any cohomology class $\zeta\in H^2(M, {\Bbb R})$,
we thus have a $\gamma$-induced decomposition
$$\zeta= \zeta^+ + \zeta^-, $$
where $\zeta^\pm \in {\cal H}^\pm$.

The subspace ${\cal H}^+\subset H^2(M, {\Bbb R})$
is called the {\em polarization} determined by 
$\gamma$. The intersection form $\cup : H^2\times H^2 \to H^4= {\Bbb R}$
becomes positive-definite when restricted to 
${\cal H}^+$, and ${\cal H}^+$ is a maximal subspace with
this property. Indeed,  ${\cal H}^-$ is the orthogonal 
complement of ${\cal H}^+$ with respect to 
$\cup$, and the restriction of 
$\cup$ to ${\cal H}^-$ is {\em negative}-definite. 
The dimension of ${\cal H}^\pm$ is therefore a
homotopy invariant $b^\pm$ of $M$, 
and the difference $\tau = b^+-b^-$   is
called the {\em signature} of $M$. 
It is important to point out that 
 the polarization ${\cal H}^+\subset H^2$
really does \cite{don}
 depend on the conformal class $\gamma$ unless $b^-=0$.

If $\omega$ is a self-dual harmonic 2-form 
with respect to $\gamma$, we have 
$$\|\omega\|_2^2= [\omega ] \cup [\omega ] = [\omega ]^2,$$
since $|\omega |^2d\mu = \omega\wedge \star \omega =
\omega \wedge \omega$. Thus Corollary \ref{yup}
immediately implies

\begin{prop}\label{yip}
Let $\zeta \in H^2(M, {\Bbb R})$ be a
fixed cohomology class on a smooth compact   4-manifold $M$,
and let $\gamma$ be any conformal class
on $M$. Let
$\phi$ denote the unique $\gamma$-harmonic 2-form 
with $[\phi ] = \zeta$. Then one of the following must 
hold: 
\begin{itemize}
\item  there is a $C^\infty$  metric ${g}\in \gamma $ 
with scalar curvature ${s} > |\phi^+ |_{{g}}$;
or 
\item  
$Y_\gamma <  
\sqrt{(\zeta^+)^2}$; or 
\item  
$Y_\gamma =  \sqrt{(\zeta^+)^2}$,
and there is a Yamabe   metric $g\in \gamma$
 with $s = |\zeta^+ | \equiv \mbox{\rm const}$. 
\end{itemize}
\end{prop}
\begin{proof}
Set $\omega= \phi^+$ and apply Corollary \ref{yup}. 
\end{proof}

\section{Dirac Operators and Spin$^c$ Structures}

Let $M$ be a compact oriented 4-manifold. 
A cohomology class $\eta \in H^2(M, {\Bbb Z})$ is then called
{\em characteristic} if $\eta \equiv w_2 (M)\bmod 2$;
by a theorem of Wu, such elements always exist. 
Given such a class $\eta$, let $L$ be the Hermitian
complex line bundle with $c_1(L)=\eta$. 
This $L$ is unique up to 
isomorphism.  Moreover, given a conformal class
 $\gamma$ on $M$, the 
obstruction to the existence of a 
square-root $L^{1/2}$ of $L$ precisely
coincides with the obstruction to defining
the spin bundles ${\Bbb S}_\pm$ of $(M,\gamma)$.
Thus one may define two rank-2 Hermitian vector
bundles $V_\pm$ on $M$ such that 
$$V_\pm = {\Bbb S}_\pm\otimes L^{1/2},$$
in the formal sense that on an any spin open set of
$M$,    ${\Bbb S}_\pm$ and $L^{1/2}$ may be
defined, and there is a  canonical (but sign-ambiguous) 
isomorphism $V_\pm \to {\Bbb S}_\pm\otimes L^{1/2}$.
A choice of such bundles $V_\pm$ is called a
spin$^c$ structure. If $H_1(M,{\Bbb Z})$ contains
no elements of order 2, the spin$^c$ structures
on $M$ are in one-to-one correspondence with 
the set of characteristic elements $\eta\in H^2(M, {\Bbb Z})$.

Now fix a spin$^c$ structure on $M$, and 
choose some Hermitian connection $\theta$ on the associated 
line bundle $L\to M$. If $g$ is any metric on $M$,
 its  Levi-Civit\`a
 connection and $\theta$ together induce a
connection  
$$\nabla^\theta : \Gamma (V_+) \to \Gamma (V_+\otimes T^*M)$$
via the local identifications
$$V_+={\Bbb S}_+\otimes L^{1/2}.$$
On the other hand,
   Clifford multiplication induces a 
bundle homomorphism 
$$  V_+\otimes T^*M\stackrel{\cdot}{\to} V_-.$$
Composing these maps gives us a (twisted) Dirac
operator
$$D^\theta: \Gamma (V_+)\to \Gamma (V_-).$$
The latter  is an elliptic operator 
whose index is given by 
$$\mbox {ind } (D^\theta) = \frac{c_1^2(L)-\tau (M)}{8}.$$
When this index is positive,  
we get an estimate for the Yamabe constant
of any conformal class:

\begin{thm}\label{yap} 
Let $M$ be a smooth compact oriented 4-manifold,
and let $\eta\in H^2(M,{\Bbb Z})$ 
be a non-torsion, characteristic element  such that 
$\eta^2 > \tau (M)$. 
Let $\gamma$ be  any 
smooth conformal class on $M$.
Then 
$$Y_\gamma \leq 4\pi \sqrt{2(\eta^+)^2}. $$
Moreover,  equality occurs  iff
$M$ is  diffeomorphic to a 
rational complex surface, in such a manner
that $\eta$ becomes the first Chern class $c_1(M)$,
and  some Yamabe metric representing
$\gamma$ becomes  a  K\"ahler metric
of  constant, non-negative scalar curvature. 
\end{thm}
\begin{proof} 
Let $\varphi$ denote the  unique  $\gamma$-harmonic
2-form  such that the de Rham class
$[\varphi ]$ coincides with the image of 
 $\eta$ in real cohomology. 
If we had $Y_\gamma \geq 4\pi \sqrt{2 (\eta^+)^2}$,
Proposition \ref{yip},
applied to   $\zeta = 4\pi \sqrt{2} \eta$,
would assert the existence of 
a smooth  metric    $g\in\gamma$ with 
$s_g \geq 4\pi \sqrt{2} |\varphi^+ |_g$; and 
if equality holds, moreover, $g$ may be
further assumed to be a Yamabe metric. 
We claim, however, that this leads to 
a contradiction unless equality holds and 
the geometry is of the special kind detailed above.

Indeed, 
set $F= -2\pi i\varphi$, and let 
 $L\to M$ be the unique Hermitian line bundle
with $c_1(L)=\eta$. 
Since $\frac{i}{2\pi} F$ then represents 
the image of $c_1(L)$ in real cohomology,
the Chern-Weil theorem tells us there
is a $U(1)$ connection $\theta$ on $L$ 
whose curvature is $F$. Choose a spin$^c$ 
structure with associated line bundle $L$,
and let $D^\theta :   \Gamma (V_+)\to \Gamma (V_-)$
be the corresponding Dirac operator. 
By construction, the
index 
$$\mbox {ind } (D^\theta) = \frac{\eta^2-\tau (M)}{8}$$
of this operator 
is positive.  Thus there exists a   smooth section 
$\psi\not\equiv 0$ of $V_+$  with $D^\theta\psi =0$.
But, by the Weitzenb\"ock formula \cite{hit,LM} 
$$D^{\theta*}D^\theta
 = \nabla^{\theta *}\nabla^\theta + \frac{s}{4} + \frac{1}{2}F^+,$$
and we therefore have 
$$0= (\psi, \nabla^*\nabla \psi) + \frac{s}{4}|\psi |^2 
+ \frac{1}{2}(\psi, F^+\cdot \psi),$$
where the   self-dual 2-form
 $F^+$  acts  on $V_+$ by Clifford multiplication.
The latter action 
  is diagonalizable, with eigenvalues $\pm \sqrt{2}|F^+|=\pm
2\pi \sqrt{2}|\varphi^+ |$.
Thus
$$0\geq (\psi, \nabla^*\nabla \psi) + \frac{s- 4\pi \sqrt{2}|\varphi^+|}{4} 
|\psi |^2.$$
Integrating  over $M$, we thus have
$$0 \geq \int_M [|\nabla \psi |^2 + \frac{s- 4\pi \sqrt{2}|\varphi^+|}{4}
|\psi |^2 ]  d\mu.$$
But, by assumption, our metric $g$
satisfies  $s \geq  4\pi \sqrt{2}|\varphi^+|$. It
follows that $s =  4\pi \sqrt{2}|\varphi^+|$,
$g$ is a Yamabe metric, and  
 $\nabla\psi =0$.  The non-zero self-dual 2-form 
$\psi \odot \bar{\psi}$ is therefore parallel, 
 so $g$ is  a
K\"ahler metric. Moreover, $L$ is now the anti-canonical line bundle of 
the associated complex structure on 
$M$, so $\eta = c_1(M)$. But since 
$\eta$ is not a torsion class
by assumption, and since $M$ 
admits a K\"ahler metric of non-negative scalar curvature, 
our  complex surface
 is rational or ruled \cite{yau}.  Moreover, its
 Todd genus $(c_1^2-\tau)/8$   is also  positive, so
it follows \cite{bpv} that 
$M$ is rational --- i.e. obtained from 
${\Bbb CP}_2$ by blowing up and down. 
\end{proof}

\section{The Main Theorems} 

We now restrict the last result to 4-manifolds with 
positive-definite intersection forms.

\begin{thm}
Let $M$ be a smooth compact oriented 4-manifold
with $b_-(M) =0$,  
and suppose that  $\eta\in H^2(M,{\Bbb Z})$ 
is a characteristic element  such that 
$\eta^2 > b_2 (M)$. 
Then 
$$Y (M) \leq 4\pi \sqrt{2\eta^2}. $$
Moreover, if there is a 
a conformal class $\gamma$ on 
$M$ such that $Y_\gamma = 4\pi \sqrt{2\eta^2}$,
then $\eta^2=9$ and 
$M$ is  diffeomorphic to ${\Bbb CP}_2$
in 
such a manner
that  
$\gamma$ becomes  the conformal class of the
Fubini-Study metric. 
\end{thm}
\begin{proof} 
Because $b_-(M)=0$, $\eta^+=\eta$
for any conformal class $\gamma$, and
Theorem \ref{yap} therefore asserts that 
$Y_\gamma \leq 4\pi \sqrt{2\eta^2}$
for any conformal class $\gamma$. 
Taking the
supremum over all $\gamma$ then yields the 
first part of the result.

The second part of the result similarly follows from
Theorem \ref{yap} because ${\Bbb CP}_2$ is the
only rational surface with $b_-=0$, and
because \cite{lic0} the isometry group
of any  constant-scalar-curvature
K\"ahler metric on ${\Bbb CP}_2$ must be a
maximal compact subgroup
 of the complex automorphism group $PGL(3, {\Bbb C})$.
\end{proof}

\begin{cor} 
 Let $M$ be a smooth compact 4-manifold with 
non-trivial,  positive-definite  
intersection form.  Then $Y(M) \leq 4\pi \sqrt{2b_2(M)+16}$.
\end{cor}
\begin{proof}
We may assume that $b_2 < 4$, since the 
 the upper bound in question is otherwise a trivial 
consequence of Aubin's estimate. Thus the intersection
form of $M$ is automatically diagonalizable \cite[p.19]{HM}. 
Choose a basis for the free part of $H^2$
relative to which the intersection form 
is represented by the identity matrix,
and let  $\eta$ be a characteristic element 
whose  free part  is  a truncation of $(3,1, 1)$
in this basis. Then $\eta^2 = 8+b_2(M) > b_2(M)$.
Now apply the previous theorem. 
\end{proof}
 
\begin{thm}
Let $X_1, X_2, \ldots, X_\ell$ be    3-dimensional
spherical 
space-forms, and let 
$$M= k{\Bbb CP}_2 \# (S^1 \times X_1 ) \# \cdots
\# (S^1\times X_\ell )$$
for some $k \geq 1$. Then 
$$12\sqrt{2}\pi \leq Y(M) \leq 4\pi \sqrt{2k+16}.$$
\end{thm}
\begin{proof}
The upper bound in question is precisely that
provided by the previous 
corollary. 

 To obtain the lower bound, 
first let $g$ denote the Fubini-Study metric
on ${\Bbb CP}_2$. This is an Einstein metric, and hence a 
Yamabe minimizer by the 4-dimensional 
Gauss-Bonnet theorem. Thus 
$$Y({\Bbb CP}_2) \geq Y_{[g]} = {\cal S}(g) = 12\sqrt{2}\pi.$$
Next,  recall    \cite{okob,sch2} 
that $Y(S^1\times X_j) =Y(S^4)$. 
Now a  fundamental result 
of O. Kobayashi \cite{okob} asserts that 
$$Y(M_j)  \geq 0 ~\forall j ~ 
\Longrightarrow ~ Y(M_1\# \cdots \# M_n) \geq  \min_j Y(M_j), $$
so we therefore have 
$$Y( k{\Bbb CP}_2 \# (S^1 \times X_1 ) \# \cdots
\# (S^1\times X_\ell )) \geq Y ({\Bbb CP}_2) \geq 12\sqrt{2}\pi , $$
which is precisely the promised lower bound.  
\end{proof}

Theorems A and B are simply   interesting
special cases of this
result. 

Theorem C also follows  quite 
easily. Indeed, a primary  Hopf surface 
is  diffeomorphic to 
$S^1\times S^3$, whereas
its one-point blow-up is diffeomorphic
 to
${\Bbb CP}_2\# (S^1\times S^3)$,
albeit in
an orientation-reversing manner. 
Thus a primary Hopf surface has Yamabe
invariant equal to  $Y(S^4) = 8\sqrt{6} \pi$
by \cite{okob,sch2}, whereas its blow-up 
has Yamabe invariant equal to
$Y({\Bbb CP}_2) = 12\sqrt{2}\pi$ by the above result. 
Incidentally, the same argument  also works for secondary 
Hopf surfaces (finite quotients
of primary Hopf surfaces).

\end{document}